\documentclass[11pt,leqno,a4paper]{article}

\usepackage{texdraw,multicol,rotating}
\usepackage{amsmath}
\usepackage{epsfig}
\usepackage{texdraw}
\usepackage{amssymb,amsthm,amsfonts,amscd}
\input{txdtools.tex}
\theoremstyle{plain}
\newtheorem{thm}{Theorem}[section]
\newtheorem{prop}[thm]{Proposition}
\newtheorem{lem}[thm]{Lemma}
\newtheorem{cor}[thm]{Corollary}
\theoremstyle{definition}
\newtheorem{defi}[thm]{Definition}

\newtheorem{example}[thm]{Example}
\theoremstyle{remark}
\newtheorem{remark}[thm]{Remark}
\numberwithin{equation}{section}
\newbox{\tmpa}
\newbox{\tmpb}

\newcommand{\nc}{\newcommand}
\nc{\Uq}{U_q} \nc{\Z}{\mathbf{Z}} \nc{\C}{\mathbf{C}}
\nc{\Q}{\mathbf{Q}}\nc{\N}{\mathbf{N}}
\nc{\op}{\oplus} \nc{\ot}{\otimes} \nc{\pv}{P^{\vee}}
\nc{\ali}{\alpha_i} \nc{\B}{\mathbf{B}} \nc{\F}{\mathbf{F}}
\nc{\bP}{\mathbf{P}} \nc{\V}{\mathbf{V}} \nc{\La}{\Lambda}
\nc{\la}{\lambda}
\nc{\nbinom}[2]{\genfrac{}{}{0pt}{1}{{#1}}{{#2}}}
\nc{\qbinom}[2]{\left[\genfrac{}{}{0pt}{1}{{#1}}{{#2}}\right]}
\nc{\path}{\mathcal{P}} \nc{\fit}{\tilde{f}_i}
\nc{\eit}{\tilde{e}_i} \nc{\Y}{\mathbf{Y}} \nc{\A}{\mathbf{A}}
\nc{\ra}{\rightarrow} \nc{\vep}{\varepsilon} \nc{\vphi}{\varphi}
\nc{\g}{\mathfrak{g}} \nc{\h}{\mathfrak{h}} \nc{\oP}{\overline{P}}
\nc{\pathp}{\mathbf{p}}
\nc{\tris}{ \bsegment \move(0 0)\lvec(10 0)\lvec(10 10)\lvec(0
0)\ifill f:0.7 \esegment } \nc{\recs}{ \bsegment \move(0
0)\lvec(10 0)\lvec(10 5)\lvec(0 5)\lvec(0 0)\ifill f:0.7 \esegment
}
\nc{\hcvec}[5]{%
\getpos(#1 #3)\spx\spy \getpos(#2 #3)\epx\epy \getpos(#4
#5)\xoff\yoff \realadd \spx \xoff \twox \realadd \epx {-\xoff}
\thrx \realadd \spy \yoff \posy \move({\spx} {\spy}) \clvec
({\twox} {\posy})({\thrx} {\posy})({\epx} {\epy}) \rmove(0 0) }
\nc{\ahead}[2]{%
\cossin (0 0)({#1} {#2})\cosa\sina \bsegment
  \drawdim in \setunitscale 0.065
  \realmult {-0.5} \cosa \hcosa
  \realmult {-0.5} \sina \hsina
  \move({\hcosa} {\hsina}) \ravec({\cosa} {\sina})
\esegment }
\nc{\boxi}{%
{%
\savebox{\tmppic}{\begin{texdraw} \small \drawdim em \textref h:C
v:C \setunitscale 0.55 \htext(0 0){$i$} \move(-1 -1)\lvec(-1
1)\lvec(1 1)\lvec(1 -1)\lvec(-1 -1)
\end{texdraw}}%
\raisebox{-0.19\height}{\usebox{\tmppic}}%
}%
}
\nc{\boxj}{%
{%
\savebox{\tmppic}{\begin{texdraw} \small \drawdim em \textref h:C
v:C \setunitscale 0.55 \htext(0 0.1){$j$} \move(-1 -1)\lvec(-1
1)\lvec(1 1)\lvec(1 -1)\lvec(-1 -1)
\end{texdraw}}%
\raisebox{-0.19\height}{\usebox{\tmppic}}%
}%
}
\nc{\boxipo}{%
{%
\savebox{\tmppic}{\begin{texdraw} \small \drawdim em \textref h:C
v:C \setunitscale 0.55 \htext(0.15 0){$i\!\!+\!\!1$} \move(-1.4
-1)\lvec(-1.4 1)\lvec(1.4 1)\lvec(1.4 -1)\lvec(-1.4 -1)
\end{texdraw}}%
\raisebox{-0.19\height}{\usebox{\tmppic}}%
}%
} \everytexdraw{ \drawdim in \arrowheadsize l:0.065 w:0.03
\arrowheadtype t:F \textref h:C v:C }
\newsavebox{\tmppic}
\newsavebox{\tmpfig}
\newsavebox{\tmpdraw}
\newsavebox{\tmpfiga}
\newsavebox{\tmpfigb}
\newsavebox{\tmpfigc}
\newsavebox{\tmpfigd}
\newsavebox{\tmpfige}
\newsavebox{\tmpfigf}
\newsavebox{\tmpfigg}
\newsavebox{\tmpfigh}
\newsavebox{\tmpfigi}
\newsavebox{\tmpfigj}
\newsavebox{\tmpfigk}
\newsavebox{\tmpfigl}
\newsavebox{\tmpfigm}
\newsavebox{\tmpfign}
\newsavebox{\tmpfigo}
\newsavebox{\tmpfigp}
\newsavebox{\tmpfigq}
\newsavebox{\tmpfigr}
\newsavebox{\tmpfigs}
\newsavebox{\tmpfigt}
\newsavebox{\tmpfigu}
\newsavebox{\tmpfigv}
\newsavebox{\tmpfigw}
\newsavebox{\tmpfigx}
\newsavebox{\tmpfigy}
\newsavebox{\tmpfigz}
\newsavebox{\tmpfigaa}
\newsavebox{\tmpfigab}
\newsavebox{\tmpfigac}
\newsavebox{\tmpfigad}
\newsavebox{\tmpfigae}
\newsavebox{\tmpfigaf}
\newsavebox{\tmpfigag}
\newsavebox{\tmpfigah}
\newsavebox{\tmpfigai}
\newsavebox{\tmpfigaj}
\newsavebox{\tmpfigak}
\newsavebox{\tmpfigal}
\newsavebox{\tmpfigam}
\newsavebox{\tmpfigan}
\newsavebox{\tmpfigao}
\newsavebox{\tmpfigap}
\newsavebox{\tmpfigaq}
\newsavebox{\tmpfigar}
\newsavebox{\tmpfigas}
\newsavebox{\tmpfigat}
\newsavebox{\tmpfigau}
\newsavebox{\tmpfigav}
\newsavebox{\tmpfigaw}
\newsavebox{\tmpfigax}
\newsavebox{\tmpfigay}
\newsavebox{\tmpfigaz}
\newsavebox{\tmpfigba}
\newsavebox{\tmpfigbb}
\newsavebox{\tmpfigbc}
\newsavebox{\tmpfigbd}
\newsavebox{\tmpfigbe}
\newsavebox{\tmpfigbf}
\newsavebox{\tmpfigbg}
\newsavebox{\tmpfigbh}

\newsavebox{\spinaa}
\newsavebox{\spinab}
\newsavebox{\spinac}
\newsavebox{\spinad}
\newsavebox{\spinae}
\newsavebox{\spinaf}
\newsavebox{\spinag}
\newsavebox{\spinah}
\newsavebox{\spinai}
\newsavebox{\spinaj}
\newsavebox{\spinak}
\newsavebox{\spinal}
\newsavebox{\spinam}
\newsavebox{\spinba}
\newsavebox{\spinbb}
\newsavebox{\spinbc}
\newsavebox{\spinbd}
\newsavebox{\spinbe}
\newsavebox{\spinbf}
\newsavebox{\spinbg}
\newsavebox{\spinbh}
\nc{\node}{\lcir r:1 }
\nc{\sline}{\bsegment\savepos(10 0)(*ex *ey)
            \move(1 0)\rlvec(8 0)
            \esegment\move(*ex *ey)}
\nc{\dline}{\bsegment\savepos(10 0)(*ex *ey)
            \move(0.93 0.4)\rlvec(8.14 0)\rmove(0 -0.8)\rlvec(-8.14 0)
            \esegment\move(*ex *ey)}
\nc{\uline}{\bsegment\savepos(0 10)(*ex *ey)
            \move(0 1)\rlvec(0 8)
            \esegment\move(*ex *ey)}
\nc{\lpoint}{\savecurrpos(*ex *ey)
             \rmove(2.5 1.5)\rlvec(-1.5 -1.5)\rlvec(1.5 -1.5)
             \move(*ex *ey)}
\nc{\rpoint}{\savecurrpos(*ex *ey)
             \rmove(-2.5 -1.5)\rlvec(1.5 1.5)\rlvec(-1.5 1.5)
             \move(*ex *ey)}
\nc{\bline}{\bsegment\savepos(10 0)(*ex *ey)
            \linewd 0.6 \move(1.1 0)\rlvec(7.8 0)
            \esegment\move(*ex *ey)}
\nc{\araise}[1]{\raisebox{4.5pt}{#1}}
\nc{\braise}[1]{\raisebox{12.1pt}{#1}}
\nc{\craise}[1]{\raisebox{8pt}{#1}}
\nc{\draise}[1]{\raisebox{12pt}{#1}}
\nc{\eraise}[1]{\raisebox{14.5pt}{#1}}

\begin{document}

\title{\bf Crystal Bases and Monomials for
$U_q(G_2)$-modules}

\author{Dong-Uy Shin
\begin{thanks}
{This research was supported by KOSEF Grant \# 98-0701-01-5-L}
\end{thanks}\cr {
} \cr {\small School of Mathematics} \cr {\small Korea Institute
for Advanced Study} \cr {\small Seoul 130-722, Korea} \cr
{\small\texttt{shindong@kias.re.kr}}}
\date{}
\maketitle

\begin{abstract}
In this paper, we give a new realization of crystal bases for
irreducible highest weight modules over $U_q(G_2)$ in terms of
monomials. We also discuss the natural connection between the
monomial realization and tableau realization.
\end{abstract}
\maketitle
\vskip 1cm

\section*{Introduction}
In 1985, the {\it quantum groups} $U_q(\frak g)$, which may be
thought of as $q$-deformations of the universal enveloping
algebras $U(\frak g)$ of Kac-Moody algebras $\frak g$, were
introduced independently by Drinfel'd and Jimbo \cite{Drin,Jim}.
The integrable highest weight representations over symmetrizable
Kac-Moody algebras can be deformed consistently to the highest
weight representations over the corresponding quantum groups for
generic $q$ \cite{Lu}. From this point of view, the {\it crystal
basis theory} for integrable modules over quantum groups was
developed by Kashiwara \cite{Kas90,Kas91}. Crystal bases can be
viewed as bases at $q=0$ and they are given a structure of colored
oriented graphs, called the {\it crystal graphs}. Crystal graphs
have many nice combinatorial properties reflecting the internal
structure of integrable modules. In \cite{KasNak}, Kashiwara and
Nakashima gave an explicit realization of crystal bases of finite
dimensional irreducible modules over classical Lie algebras using
semistandard tableaux with given shapes satisfying certain
additional conditions. Motivated by their work, Kang and Misra
discovered a tableau realization of crystal bases for finite
dimensional irreducible modules over the exceptional Lie algebra
$G_{2}$ \cite{KM}. In \cite{Li1}, Littelmann gave another
description of crystal bases for finite dimensional simple Lie
algebras using the Lakshmibai-Seshadri monomial theory. His
approach was generalized to the {\it path model theory} for all
symmetrizable Kac-Moody algebras \cite{Li2,Li3}.

In \cite{KasSai}, Kashiwara and Saito gave a geometric realization
of the crystal graph $B(\infty)$ of $U_q^-(\frak g)$ as the set of
irreducible components of a lagrangian subvariety ${\mathcal L}$
of the quiver variety $\frak M$  and in \cite{Saito}, Saito
extended their idea to the crystal base $B(\la)$ of irreducible
highest weight modules of $U_q(\frak g)$. In \cite{Nakaj1}, a
crystal structure on the set of irreducible components of a
lagrangian subvariety $\frak Z$ of the quiver variety $\frak M$
was given. In \cite{Nakaj2}, Nakajima gave the set of $U_q({\bf
L}\frak g)$-modules $M(P)$, called the {\it standard modules},
with a loop algebra ${\bf L}\frak g$ of $\frak g$ and Drinfel'd
polynomials $P$, and in \cite{Nakaj3} he introduced the
$t$-analogs $\chi_{q,t}$ of $q$-character. While studying the
$t$-analogs of $q$-character of standard modules $M$, he
discovered that these irreducible components of a lagrangian
subvariety $\frak Z$ are identified with certain monomials, and so
the action of Kashiwara operators can be interpreted as
multiplication by monomials. Moreover, in \cite{Kas02,Nakaj4},
Kashiwara and Nakajima gave a crystal structure on the set
${\mathcal M}$ of monomials and they showed that the connected
component ${\mathcal M}(\lambda)$ of ${\mathcal M}$ containing a
highest weight vector $M$ with a dominant integral weight
$\lambda$ is isomorphic to the irreducible highest weight crystal
$B(\la)$. But, they did not give an explicit characterization of
monomials in ${\mathcal M}(\lambda)$. For the special linear Lie
algebra, in \cite{KKS}, Kang, Kim and the Author gave an explicit
characterization and investigated the connection with the
realization in terms of semistandard tableaux. Moreover, for the
affine Lie algebra $A_n^{(1)}$, in \cite{Kim}, Kim gave an
explicit description and she gave crystal isomorphisms from
monomial realization to path realization and Young wall
realization given in \cite{HKL,JMMO,Kang,KMN1}.

In this paper, for any dominant integral weight $\lambda$, we give
an explicit description of the crystal ${\mathcal M}(\lambda)$ for
$U_q(G_2)$. In addition, we discuss the connection between the
monomial realization and tableau realization of crystal bases
given by Kang and Misra.

\vskip 3mm {\bf Acknowledgments.} \ The author would like to
express his sincere gratitude to Professor S.-J. Kang and Doctor
J.-A. Kim for their interest in this work and many valuable
discussions.


\vskip 1cm
\section{Crystal bases and Nakajima's monomial}
\subsection{Crystal bases for $U_q(G_2)$}
Let $A=\left(%
\begin{array}{rr}
  2 & -3 \\
  -1 & 2 \\
\end{array}%
\right)$ be a Cartan matrix of type $G_2$. Then we have the finite
dimensional simple Lie algebra $G_2$ and the {\it quantum group}
$U_q(G_2)$ associated to $A$ with a {\it Cartan subalgebra} $\frak
h$, the set of {\it simple roots} $\Pi=\{\alpha_1,\alpha_2\}$ and
the set of {\it simple coroots} $\Pi^{\vee} = \{ h_1, h_2\}$. Let
$\Lambda_i\in \frak h^*$ ($i=1,2$) be the {\it fundamental weight}
and set $P=\Z \Lambda_1\oplus\Z \Lambda_2$ and $P^+=\Z_{\ge 0}
\Lambda_1\oplus\Z_{\ge 0} \Lambda_2$. Then every finite
dimensional $U_q(G_2)$-module $M$ is a direct sum of irreducible
modules $V(\la)$ with $\la\in P^+$.

Let $M$ be a finite dimensional $U_q(\frak g)$-module, then we
have the {\it crystal basis} $(L, B)$ \cite{Kas90,Kas91}. The set
$B$ gives a colored oriented graph structure with the arrow
defined by
\begin{center}
$b \stackrel{i} \longrightarrow b' \quad \text{if and only if}
\quad \fit b = b'.$
\end{center}
The graph $B$ is called the {\it crystal  graph} of $M$. Moreover,
the crystal bases have a nice behavior with respect to the tensor
product. See \cite{HK} for more details of quantum groups and
crystal bases.

 Now we give the description of
the set of tableaux of type $G_2$ which is a realization of
crystal basis \cite{KM}. Let $\la$ be a dominant integral weight.
Then the crystal graph $B(\la)$ of the irreducible highest weight
module $V(\la)$ is realized as the set of tableaux on
$I=\{1,2,3,0,\overline{3},\overline{2},\overline{1}\}$ with the
linear order $1\prec 2\prec 3\prec 0\prec \overline{3}\prec
\overline{2}\prec \overline{1}.$ To describe $B(\la)$, we need
some definitions and conditions.

\begin{defi}
For $i,j=0,1,2,3$, we define
$$
\begin{array}{ll}
\text{dist}(i,j)=\text{dist}(\overline{j},\overline{i})=j-i
&\text{for $i<j,\,\,i,j=1,2,3$},\\
\text{dist}(i,0)=\text{dist}(0,\overline{i})=4-i
&\text{for $i=1,2,3$},\\
\text{dist}(i,\overline{j})=8-(i+j)
&\text{for $i,j=1,2,3$}.\\
\end{array}
$$
We also define $\text{dist}(a,b)=\text{dist}(b,a)$ and
$\text{dist}(a,a)=0$ for all $a,b\in I$.
\end{defi}

Now, we are ready to give a characterization of tableaux for type
$U_q(G_2)$.

\begin{prop} For a dominant integral weight $\la$, the crystal graph $B(\la)$
is realized as the set of tableaux $T$ of shape $\la$ with entries
on $I$ such that
\begin{itemize}
\item [{\rm (i)}] the entries of $T$ weakly increase along the
rows, but the element $0$ cannot appear more than once,

\item [{\rm (ii)}] the entries of $T$ strictly increase down the
columns, but the element $0$ can appear more than once,

\item [{\rm (iii)}] for each column $C$ of length $2$ with the
entries $a$,$b$, $\text{dist}(a,b)\le 2$ for $a=1,0$ and
$\text{dist}(a,b)\le 3$ otherwise,

\item [{\rm (iv)}] for each pair of adjacent columns $C$, $C'$ of
length $2$ with the entries (from left to right and from top to
bottom) $a,b,c,d$, $\text{dist}(a,d)\ge 3$ for $a=2,3,0$ and
$\text{dist}(a,d)\ge 2$ for $a=\overline{3}$.
\end{itemize}
\end{prop}

\vskip 5mm
\subsection{Nakajima's monomials} In this subsection, we
recall the crystal structure on the set of monomials discovered by
Nakajima \cite{Nakaj4}. Our presentation follows that of Kashiwara
\cite{Kas02}. Moreover,  we only treat the  monomials for the type
$U_q(G_2)$.

Let $\mathcal M$ be the set of monomials in the variables $Y_i(n)$
for $i=1,2$ and $n\in {\bf Z}$. Here, a typical elements $M$ of
${\mathcal M}$ has the form $M=Y_{i_1}(n_1)^{a_1}\cdots
Y_{i_r}(n_r)^{a_r}$, where $i_k=1,2$ and $n_k, a_k \in {\bf Z}$
for $k=1,\cdots, r.$ Since $Y_i(n)$'s are commuting variables, we
may assume that $n_1\leq n_2\leq\cdots \leq n_r.$

For a monomial $M=Y_{i_1}(n_1)^{a_1}\cdots Y_{i_r}(n_r)^{a_r}$, we
define
\begin{center}
$ \aligned
\text{wt}(M)&=\sum_{k=1}^{r}a_k\La_{i_k}=a_1\La_{i_1}+\cdots a_r\La_{i_r},\\
\varphi_i(M)&=\text{max}\,\big(\big\{\sum_{k=1 \atop i_k=i}^s a_k \,\,|\,\,1\leq s\leq r \big\}\cup \{0\}\big),\\
\varepsilon_i(M)&=\text{max}\,\big(\big\{-\sum_{k=s+1 \atop
i_k=i}^r a_k \,\,|\,\,1\leq s\leq r-1 \big\}\cup \{0\}\big).
\endaligned
$
\end{center}
It is easy to verify that $\varphi_i(M)\geq 0,\varepsilon_i(M)\geq
0,$ and $\langle h_i,
\text{wt}(M)\rangle=\varphi_i(M)-\varepsilon_i(M)$.

First, we define
$$
\aligned
n_f&=\text{smallest}\,\,n_s\,\, \text{such that}\,\, \varphi_i(M)=\sum_{k=1 \atop i_k=i}^s a_k\\
n_e&=\text{largest}\,\,n_s\,\, \text{such that}\,\, \varphi_i(M)=\sum_{k=1 \atop i_k=i}^s a_k.\\
\endaligned
$$
In addition, choose integers $c_{12}$ and $c_{21}$ such that
$c_{12}+c_{21}=1,$ and define
$$
\aligned A_1(n)&=Y_1(n)Y_1(n+1)Y_2(n+c_{21})^{-1},\\
A_2(n)&=Y_2(n)Y_2(n+1)Y_1(n+c_{12})^{-3}.\endaligned$$

Now, the {\it Kashiwara operators} $\eit$, $\fit$ ($i=1,2$) on
${\mathcal M}$ are defined as follows.
$$
\aligned &\fit(M)=
  \begin{cases}
    0 & \text{if $\varphi_i(M)=0$}, \\
    A_i(n_f)^{-1}M & \text{if $\varphi_i(M)>0$},
  \end{cases}\\
&\eit(M)=
  \begin{cases}
    0 & \text{if $\varepsilon_i(M)=0$}, \\
    A_i(n_e)M & \text{if $\varepsilon_i(M)>0$}.
  \end{cases}
\endaligned
$$
Then the maps $\text{wt}: {\mathcal M}\rightarrow P$, $\varphi_i,
\varepsilon_i:{\mathcal M}\rightarrow \Z\cup\{-\infty\},$
$\eit,\fit:{\mathcal M}\rightarrow {\mathcal M}\cup\{0\}$ define a
$U_q(\mathfrak g)$-crystal structure on ${\mathcal M}$
\cite{Kas02}.

Moreover, we have

\begin{prop}\, {\rm \cite{Kas02}} \,\,\label{prop:mono}

{\rm (a)} For each $i=1,2$, ${\mathcal M}$ is isomorphic to a
crystal graph of an integrable $U_{(i)}$-module.

{\rm (b)} Let $M$ be a monomial of weight $\la$ such that $\eit
M=0$ for all $i=1,2$, and let ${\mathcal M}(\la)$ be the connected
component of ${\mathcal M}$ containing $M$. Then there exists a
crystal isomorphism  $${\mathcal M}(\la) {\stackrel
{\sim}{\longrightarrow}} B(\la) \,\,\,\text{given by}\,\,\,
M\longmapsto v_{\la}.$$
\end{prop}

\begin{example}
Let $c_{12}=1$ and $c_{21}=0$. Then the connected component
containing $Y_1(0)$ which is isomorphic to the crystal graph
${\mathcal M}(\La_1)$ is given as follows:
\begin{center}
$\aligned Y_1(0)\stackrel{1}{\longrightarrow} Y_1(1)^{-1}Y_2(0)
&\stackrel{2}{\longrightarrow} Y_1(1)^2Y_2(1)^{-1}
\stackrel{1}{\longrightarrow} Y_1(1)Y_1(2)^{-1}\\
&\stackrel{1}{\longrightarrow} Y_1(2)^{-2}Y_2(1)
\stackrel{2}{\longrightarrow} Y_1(2)Y_2(2)^{-1}
\stackrel{1}{\longrightarrow} Y_1(3)^{-1}\endaligned$
\end{center}

\end{example}


\vskip 1cm
\section{Characterization of ${\mathcal M}(\la)$ and connection with tableau realization}
\subsection{Characterization of ${\mathcal M}(\la)$}
In this subsection, we give an explicit characterization of the
crystal ${\mathcal M}(\la)$ for $U_q(G_2)$. For simplicity, we
take integers $c_{12}=1$ and $c_{21}=0$. Then we have
$$A_1(n)=Y_1(n)Y_1(n+1)Y_2(n)^{-1},\quad A_2(n)=Y_2(n)Y_2(n+1)Y_1(n+1)^{-3}.$$

To characterize ${\mathcal M}(\la)$, we first focus on the case
that $\la=\La_k$ ($k=1,2$). Let $M_0=Y_k(m)$ for $m\in \Z$, then
we see that
$$\text{wt}(M_0)=\La_k,\,\, \varphi_i(M_0)=\delta_{ik}\,\,\text{and}\,\,
\varepsilon_i(M_0)=0\,\,\,\text{for all}\,\,i=1,2.$$ Hence $\eit
M_0=0$ for all $i=1,2$ and the connected component containing
$M_0$ is isomorphic to $B(\La_k)$ over $U_q(G_2).$ For simplicity,
we will take $M_0=Y_k(1)$, and that does not make much difference.

We set $Y_0(m)^{\pm 1}=1$ for all $m\in \Z.$ For $m\in \Z$, we
introduce new variables
\begin{equation}
\aligned
X_i(m)&=\left\{\begin{array}{ll}Y_{i-1}(m+1)^{-1}Y_i(m) &\text{for $i=1,2$,}\\
Y_1(m+1)^2Y_2(m+1)^{-1} &\text{for $i=3$,}
\end{array}\right.\\
X_{{\overline i}}(m)&=\left\{\begin{array}{ll}
Y_{i-1}(m+(4-i))Y_i(m+(4-i))^{-1} &\text{for
$i=1,2$,}\\
Y_1(m+2)^{-2}Y_2(m+1) &\text{for $i=3$,}
\end{array}\right.\\
X_0(m)&=Y_1(m+1)Y_1(m+2)^{-1}.
\endaligned
\end{equation}

Now, we give a lemma which plays a prominent role in
characterizing the connected component containing a monomial $M$.
By direct calculation, we have

\begin{lem}\label{lem:1} For $m\in \Z$, we
have

{\rm (a)}\, $X_1(m)X_0(m-1)=X_2(m)X_3(m-1)$ \quad{\rm (b)}\,
$X_1(m)X_{\overline 3}(m-1)=X_2(m)X_0(m-1)$

{\rm (c)}\, $X_1(m)X_{\overline 2}(m-1)=X_3(m)X_0(m-1)$ \quad{\rm
(d)}\, $X_1(m)X_{\overline 1}(m-1)=X_0(m)X_0(m-1)$

{\rm (e)}\, $X_2(m)X_{\overline 2}(m-1)=X_3(m)X_{\overline
3}(m-1)$ \quad{\rm (f)}\, $X_2(m)X_{\overline
1}(m-1)=X_0(m)X_{\overline 3}(m-1)$

{\rm (g)}\, $X_3(m)X_{\overline 1}(m-1)=X_0(m)X_{\overline
2}(m-1)$ \quad{\rm (h)}\, $X_0(m)X_{\overline 1}(m-1)=X_{\overline
3}(m)X_{\overline 2}(m-1)$
\end{lem}

\begin{remark}
If we consider the monomials $X_{\alpha}(m)X_{\beta}(m-1)$ in the
left hand side, we have
\begin{center}
$\begin{array}{ll} {\rm dist}(\alpha,\beta)\ge 3 &\text{if\,
$\alpha=2,3$ or $0$,}\\
{\rm dist}(\alpha,\beta)\ge 2 &\text{if\, $\alpha={\overline 3}$.}
\end{array}$
\end{center}
Similarly, if we consider the monomials
$X_{\gamma}(m)X_{\delta}(m-1)$ in the right hand side, we have
\begin{center}
$\begin{array}{ll} {\rm dist}(\gamma,\delta)\le 2 &\text{if\,
$\gamma=1,0$,}\\
{\rm dist}(\gamma,\delta)\le 3 &\text{otherwise.}
\end{array}$
\end{center}
\end{remark}

\begin{prop} \label{prop:one-column}
{\rm (a)} Let $M_0=Y_1(1)=X_1(1)$ be a monomial of weight $\La_1$
such that $\eit M_0=0$ for $i=1,2$, then  the connected component
${\mathcal M}(\La_1)$ of ${\mathcal M}$ containing $M_0$ is
$$\{X_{a}(1) \mid a=1,2,3,0,\overline{3},\overline{2},\overline{1}\}.$$

{\rm (b)} Let $M_0=Y_2(1)=X_1(2)X_2(1)$ be a monomial of weight
$\La_2$ such that $\eit M_0=0$ for $i=1,2$. Then the connected
component ${\mathcal M}(\La_2)$ of ${\mathcal M}$ containing $M_0$
is
$$\{X_{a}(2)X_{b}(1) \mid
a\prec b,\,\,\text{or $a=b=0$}\}.$$
\end{prop}

\begin{proof} Since it is very easy to see (a), we only treat (b).
By Proposition \ref{prop:mono}, it suffices to prove that ${\mathcal
M}(\La_2)\cup \{0\}$ is closed under the actions of $\eit$ and
$\fit$ for $i=1,2$, and $M_0$ is the unique monomial of weight $\La_2$
such that $\eit M_0=0$ for $i=1,2$. But, it is clear by
the following connected component of $Y_2(1)=X_1(2)X_2(1)$.

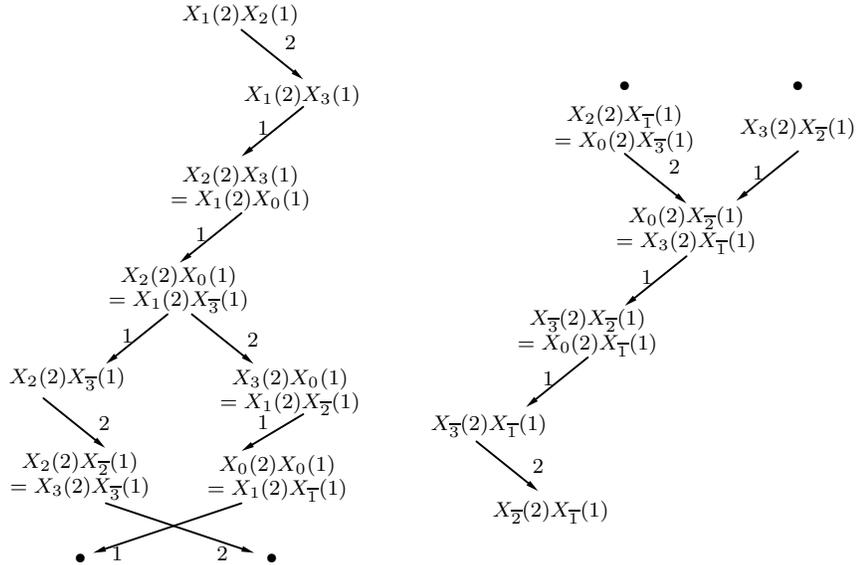
\savebox{\tmpfiga}{\begin{texdraw} \fontsize{8}{8}\selectfont
\drawdim em \setunitscale 0.7 \textref h:C v:C \setunitscale 0.55

\htext(0 0){$X_1(2)X_2(1)$}

\htext(5 -6.5){$X_1(2)X_3(1)$}

\htext(0 -13){$X_2(2)X_3(1)$} \htext(0 -15){$=X_1(2)X_0(1)$}

\htext(-5 -21.2){$X_2(2)X_0(1)$} \htext(-5
-23.2){$=X_1(2)X_{\overline 3}(1)$}

\htext(-14 -29.5){$X_2(2)X_{\overline 3}(1)$}

\htext(4 -29.5){$X_3(2)X_0(1)$} \htext(4
-31.5){$=X_1(2)X_{\overline 2}(1)$}

\htext(-13 -36.3){$X_2(2)X_{\overline 2}(1)$} \htext(-13
-38.3){$=X_3(2)X_{\overline 3}(1)$}

\htext(3 -36.5){$X_0(2)X_0(1)$} \htext(3
-38.5){$=X_1(2)X_{\overline 1}(1)$}

\htext(2.5 -44){$\bullet$}\htext(-13 -44){$\bullet$}

\move(0 -1.2) \ravec(5 -4)\rmove(-1 3)\htext{$2$}

\move(5 -7.4)\ravec(-5 -4)\rmove(1.8 2.3)\htext{$1$}

\move(0 -16)\ravec(-5 -4)\rmove(1.8 2.3)\htext{$1$}

\move(-6 -24.2)\ravec(-5 -4)\rmove(1.8 2.3)\htext{$1$}

\move(-4 -24.2) \ravec(5 -4)\rmove(0 2)\htext{$2$}

\move(-16 -31) \ravec(5 -4)\rmove(0 2)\htext{$2$}

\move(5 -32.3)\ravec(-5 -3)\rmove(1.8 2.3)\htext{$1$}

\move(0 -39.5)\ravec(-12 -4)\rmove(2 0)\htext{$1$}

\move(-11 -39.5) \ravec(12 -4)\rmove(-2.5 0)\htext{$2$}

\end{texdraw}}
\savebox{\tmpfigb}{\begin{texdraw} \fontsize{8}{8}\selectfont
\drawdim em \setunitscale 0.7 \textref h:C v:C \setunitscale 0.55

\htext(0 2.5){$\bullet$}\htext(14 2.5){$\bullet$}

\htext(0 0){$X_2(2)X_{\overline 1}(1)$}\htext(0
-2){$=X_0(2)X_{\overline 3}(1)$}

\htext(14 -1){$X_3(2)X_{\overline 2}(1)$}

\htext(5 -8.2){$X_0(2)X_{\overline 2}(1)$} \htext(5
-10.2){$=X_3(2)X_{\overline 1}(1)$}

\htext(-3 -16.5){$X_{\overline 3}(2)X_{\overline 2}(1)$} \htext(-3
-18.5){$=X_0(2)X_{\overline 1}(1)$}

\htext(-11 -25){$X_{\overline 3}(2)X_{\overline 1}(1)$} \htext(-6
-32){$X_{\overline 2}(2)X_{\overline 1}(1)$}

\move(0 -3) \ravec(5 -4)\rmove(-1 3)\htext{$2$}

\move(14 -2.8)\ravec(-5 -4)\rmove(1.8 2.3)\htext{$1$}

\move(5 -11.3)\ravec(-5 -4)\rmove(1.8 2.3)\htext{$1$}

\move(-3 -19.5)\ravec(-5 -4)\rmove(1.8 2.3)\htext{$1$}

\move(-12 -26.3) \ravec(5 -4)\rmove(0 2)\htext{$2$}

\end{texdraw}}

\vskip 3mm

$$\raisebox{-0.4\height}{\usebox{\tmpfiga}}\qquad\qquad\qquad
\raisebox{-0.4\height}{\usebox{\tmpfigb}}$$

\end{proof}

\begin{remark}
 If we take $M_0=Y_k(N)$ ($k=1,2$), then we need to replace $X_i(m)$ by
$X_i(m+N-1).$
\end{remark}

Applying Lemma \ref{lem:1}, i.e., replacing monomials in the left
hand side with the monomials in the right hand side, we have

\begin{thm}\label{thm:one-column}
Let $M_0=Y_2(1)=X_1(2)X_2(1)$ be a monomial of weight $\La_2$
such that $\eit M_0=0$ for $i=1,2$. Then the connected component
${\mathcal M}(\La_2)$ of ${\mathcal M}$ containing $M_0$ is
characterized as
\begin{center}
${\mathcal M}(\La_2)= \Biggl\{X_{a}(2)X_{b}(1)\,\bigg|
\begin{array}{l}
{\rm(i)}\,a\prec b,\,\,\text{or $a=b=0$,}\\
{\rm (ii)}\, \text{${\rm dist}(a,b)\le 2$ for $a=1,0$, ${\rm
dist}(a,b)\le 3$ otherwise.}
\end{array}
\Biggr\}.$
\end{center}
\end{thm}

We now consider the general case. By direct calculation, we have

\begin{lem}\label{lem:2} For $m\in \Z$, we
have
\begin{center}
$X_0(m)^2=X_3(m)X_{\overline 3}(m).$
\end{center}
\end{lem}

\begin{prop}
\label{prop:general} Let $\la=m\La_1+n\La_2$. Then the connected
component ${\mathcal M}(\la)$ containing
\begin{center}
$M_0=Y_1(1)^{m}Y_2(1)^{n}=X_1(1)^{m}(X_1(2)X_2(1))^{n}$
\end{center}
is expressed as the set of monomials
$$\begin{array}{c} M=X_{a_1}(2)\cdots X_{a_n}(2)X_{b_1}(1)\cdots X_{b_{m+n}}(1)
\\\quad\,\,\,\,\,
{\rm(}a_1\succeq \cdots \succeq a_n,\,\,b_1\succeq \cdots \succeq
 b_{m+n}{\rm)}
\end{array}$$  satisfying following
conditions:
\begin{equation} \label{eq:prop:general}
a_j\prec b_j \quad\text{unless $a_j=b_j=0$ for
$j=1,\cdots,n$}.\end{equation}
\end{prop}

\begin{proof} As in Proposition \ref{prop:one-column}, it suffices to prove the
following statements:

(a) For all $i=1,2$, we have $\tilde e_i {\mathcal M(\la)} \subset
{\mathcal M}(\la) \cup \{0\}$ and $\tilde f_i {\mathcal M}(\la)
\subset {\mathcal M}(\la) \cup \{0\}.$

(b) If $M\in {\mathcal M(\la)}$ and $\eit M=0$ for all $i=1,2$,
then $M=M_0.$

\vskip 2mm We first prove the statement $(a)$. Let $M$ be a
monomial of ${\mathcal M}(\la)$. Assume that $\tilde f_i M\neq 0$
for $i=1,2$. Then $\fit M$ is obtained by multiplying
$A_i(k)^{-1}$ from $M$ for some $k$. Note that $A_1(k)^{-1}$ and
$A_2(k)^{-1}$ are expressed as
\begin{equation}\label{eq:A1}\aligned
A_1(k)^{-1}&=X_1(k)^{-1}X_2(k)=X_3(k-1)^{-1}X_0(k-1)\\
&=X_0(k-1)^{-1}X_{\overline 3}(k-1)=X_{\overline
2}(k-2)^{-1}X_{\overline 1}(k-2)\endaligned\end{equation} and
\begin{equation}\label{eq:A2} A_2(k)^{-1}=X_2(k)^{-1}X_3(k)=X_{\overline
3}(k-1)^{-1}X_{\overline 2}(k-1).\end{equation}

Suppose that $\fit M$ does not satisfy the condition
(\ref{eq:prop:general}). It means that
\begin{center}
$A_i(k)^{-1}=X_p(2)^{-1}X_{p+1}(2)$ for some $p$ determined by
(\ref{eq:A1}) and (\ref{eq:A2})
\end{center} and there is no
element among $b_1,\cdots,b_n$ larger than $p+1$, where $p+1$
denote by the next element of $p$ under the ordering in $I$. But,
by the condition (\ref{eq:prop:general}), there should be a
$X_{p+1}(1)$ in $M$ and so it is a contradiction.

Similarly, we can prove that $\tilde e_i {\mathcal M(\la)} \subset
{\mathcal M}(\la) \cup \{0\}.$

\vskip 2mm To prove (b), suppose $M\in {\mathcal M(\la)}$ and
$\tilde e_i M=0$ for all $i=1,2.$ Then we can easily check that
$b_{n+1},\cdots,b_{m+n},a_1,\cdots,a_n$ should be $1$. Moreover,
$b_1,\cdots,b_n$ should be $2$. Therefore, we have
$M=X_1(1)^{m}(X_1(2)X_2(1))^{n}=Y_1(1)^{m}Y_2(1)^{n}.$
\end{proof}

\begin{remark}
For dominant integral weights $\la,$ $\mu,$ and $\tau$ such that
$\la=\mu+\tau$, we have
$${\mathcal M}(\la)=\{M_{1}M_{2}\mid M_{1}\in
{\mathcal M}(\mu),\,\,\, M_2\in {\mathcal M}(\tau) \}.$$
\end{remark}

Now, we consider a monomial $M$ in ${\mathcal M}(\la)$ with a
dominant integral weight $\la=m\La_1+n\La_2$. At first, for any
expression of
$$M=X_{a_1}(2)\cdots X_{a_n}(2)X_{b_1}(1)\cdots X_{b_{m+n}}(1),$$
we put $X_i(j)$ in the
$j$-th row from bottom so that $X_b(j)$ exists to the right hand
side of $X_a(j)$ for $a\preceq b$. That is,
$$\begin{array}
{r} X_{a_n}(2)\cdots X_{a_1}(2)\\
X_{b_{n+m}}(1)\cdots X_{b_{n+1}}(1)X_{b_n}(1)\cdots X_{b_1}(1)
\end{array}$$

Secondly, we apply the following rules:
\begin{itemize} \item [(al-1)] If there is a pair
$(X_{\alpha}(2),X_{\beta}(1))$ such that
$X_{\alpha}(2)X_{\beta}(1)$ is one of monomials in the left hand
side of Lemma \ref{lem:1}, and $X_{\alpha}(2)$ and $X_{\beta}(1)$
lie in the same column or $X_{\beta}(1)$ lies in the left hand
side of $X_{\alpha}(2)$, then we replace
$(X_{\alpha}(2),X_{\beta}(1))$ with
$(X_{\gamma}(2),X_{\delta}(1))$ which is the pair of monomials
corresponding to $X_{\alpha}(2),X_{\beta}(1)$. Moreover, if there
are several such pairs, then we apply above rule from the pair
$(X_{\alpha}(2),X_{\beta}(1))$ with the largest distance between
$X_{\alpha}(2)$ and $X_{\beta}(1)$ to the pair
$(X_{\alpha'}(2),X_{\beta'}(1))$ with the smallest distance
between $X_{\alpha'}(2)$ and $X_{\beta'}(1)$.

\item [(al-2)] If there is a pair $(X_{\gamma}(2),X_{\delta}(1))$
such that $X_{\gamma}(2)X_{\delta}(1)$ is one of monomials in the
right hand side of Lemma \ref{lem:1}, and $X_{\delta}(1)$ lies in
the right hand side of $X_{\gamma}(2)$, then we replace
$(X_{\gamma}(2),X_{\delta}(1))$ with
$(X_{\alpha}(2),X_{\beta}(1))$ which is the pair of monomials
corresponding to $X_{\gamma}(2),X_{\delta}(1)$. Moreover, if there
are several such pairs, then we apply above rule from the pair
$(X_{\gamma}(2),X_{\delta}(1))$ with the largest distance between
$X_{\gamma}(2)$ and $X_{\delta}(1)$ to the pair
$(X_{\gamma'}(2),X_{\delta'}(1))$ with the smallest distance
between $X_{\gamma'}(2)$ and $X_{\delta'}(1)$.

\item [(al-3)] If there is a pair $(X_0(p),X_{0}(p))$ ($p=1,2$),
then we replace $(X_0(p),X_{0}(p))$ by $(X_3(p),X_{\overline
3}(p))$.

\end{itemize}

From now on, we denote by $[M]$ the expression of monomial $M$
obtained by (al-1),(al-2) and (al-3).

\begin{example} \label{example:algo}(a) Let $\la$ be a dominant integral weight $\La_1+3\La_2$ and let $M$
be a monomial $Y_1(2)^3Y_1(3)^{-3}Y_1(4)^{-1}Y_2(2)^2Y_2(3)^{-1}$,
then it can be expressed as
$$
M=X_2(2)X_1(2)^2X_{\overline{1}}(1)X_{\overline{2}}(1)X_{\overline{3}}(1)X_{0}(1)
$$
and so it is a monomial of ${\mathcal M}(\la)$. Moreover, we have
$$\aligned\begin{array}{r}
X_1(2)X_1(2)X_2(2)\\
X_0(1)X_{\overline{3}}(1)X_{\overline
{2}}(1)X_{\overline{1}}(1)\end{array}&{= \atop \text{(al-1)}}
\begin{array}{r}
X_1(2)X_2(2)X_2(2)\\
X_3(1)X_{\overline{3}}(1)X_{\overline
{2}}(1)X_{\overline{1}}(1)\end{array}\\
& {= \atop \text{(al-1)}}\begin{array}{r}
X_1(2)X_2(2)X_3(2)\\
X_3(1)X_{\overline{3}}(1)X_{\overline{3}}(1)X_{\overline{1}}(1)\end{array}\\
&{= \atop \text{(al-1)}}\begin{array}{r}
X_2(2)X_2(2)X_0(2)\\
X_3(1)X_{0}(1)X_{\overline{3}}(1)X_{\overline{2}}(1)\end{array}.
\endaligned
$$
Therefore,
$$[M]=X_0(2)X_2(2)^2X_{\overline 2}(1)X_{\overline 3}(1)X_0(1)X_{3}(1).$$

(b) Let $\la$ be a dominant integral weight $2\La_2$ and let $M$
be a monomial $Y_1(2)^2Y_1(4)^{-2}$, then it can be expressed as
$$
M=X_0(2)^2X_0(1)^2$$ and so it is a monomial of ${\mathcal
M}(\la)$. Moreover, we have
$$\begin{array}{r}
X_0(2)X_0(2)\\
X_0(1)X_0(1)\end{array}{= \atop \text{(al-2)}}
\begin{array}{r}
X_1(2)X_0(2)\\
X_0(1)X_{\overline{1}}(1)\end{array}{= \atop \text{(al-1)}}
\begin{array}{r}
X_2(2)X_{\overline 3}(2)\\
X_3(1)X_{\overline{2}}(1)\end{array}.
$$
Therefore,
$$[M]=X_{\overline 3}(2)X_2(2)X_{\overline{2}}(1)X_3(1).$$

\end{example}

By the algorithm (al-1),(al-2) and (al-3), we have
\begin{thm} \label{thm:general}
Let $\la=m\La_1+n\La_2$. Then the connected component ${\mathcal
M}(\la)$ containing
\begin{center}
$M_0=Y_1(1)^{m}Y_2(1)^{n}=X_1(1)^{m}(X_1(2)X_2(1))^{n}$
\end{center}
is expressed as the set of monomials
$$\begin{array}{c} M=X_{a_1}(2)\cdots X_{a_n}(2)X_{b_1}(1)\cdots X_{b_{m+n}}(1)
\\\quad\,\,\,\,\,
{\rm(}a_1\succeq \cdots \succeq a_n,\,\,b_1\succeq \cdots \succeq
 b_{m+n}{\rm)}
\end{array}$$ satisfying following conditions:
\begin{itemize}

\item [{\rm (i)}] $0$'s can not be repeated in $a_1,\cdots,a_n$
and $b_1,\cdots,b_{m+n}\frac{}{}$, respectively, \item [{\rm
(ii)}] $a_j\prec b_j$ unless $a_j=b_j=0$ for $j=1,\cdots,n$, \item
[{\rm (iii)}] for $j=1,\cdots,n$,
\begin{center}
$\begin{array}{ll} {\rm dist}(a_j,b_j)\le 2 &\text{if\,
$a_j=1$ or $0$,}\\
{\rm dist}(a_j,b_j)\le 3 &\text{otherwise,}
\end{array}$
\end{center}
\item [{\rm (iv)}] for $j=2,\cdots,n$,
\begin{center}
$\begin{array}{ll} {\rm dist}(a_{j},b_{j-1})\ge 3 &\text{if\,
$a_j=2,3$ or $0$,}\\
{\rm dist}(a_{j},b_{j-1})\ge 2 &\text{if\, $a_j={\overline 3}$.}
\end{array}$
\end{center}
\end{itemize}
\end{thm}


\vskip 5mm \subsection{Connection with other realizations of
crystal graphs} In this subsection, we give a natural 1-1
correspondence between monomial realization and tableau
realization of crystal bases of type $U_q(G_2)$.

Let $\la=m\La_1+n\La_2$ be a dominant integral weight. The tableau
realization $T(\la)$ of $B(\la)$ given by Kang and Misra was
obtained by imbedding $B(\la)$ to $B(\La_1)^{\otimes m}\otimes
B(\La_2)^{\otimes n}$. Similarly, if we imbed $B(\la)$ to
$B(\La_2)^{\otimes n}\otimes B(\La_1)^{\otimes m}$, we can have
another tableau realization $S(\la)$. That is, for a dominant
integral weight $\la=m\La_1+n\La_2$, we can associate a diagram
which is a collection of $n$ boxes in top row and $m+n$ boxes in
bottom row in right justified rows and $B(\la)$ is realized as the
set $S(\la)$ of tableaux satisfying the same condition as $T(\la)$
on this diagram.

\begin{thm}\label{thm:connection}
Let $\la=m\La_1+n\La_2$ be a dominant integral weight. Then there
is a crystal isomorphism $\psi:{\mathcal M}(\la)\rightarrow
S(\la)$.
\end{thm}

\savebox{\tmpfiga}{\begin{texdraw} \fontsize{7}{7}\selectfont
\drawdim em \setunitscale 0.8

\rlvec(15 0) \move(0 2)\rlvec(15 0)\move(8 4)\rlvec(7 0) \move(0
0)\rlvec(0 2)\move(2.5 0)\rlvec(0 2)\move(5.5 0)\rlvec(0 2)
\move(8 0)\rlvec(0 4)\move(10 0)\rlvec(0 4)\move(13 0)\rlvec(0 4)
\move(15 0)\rlvec(0 4)

\htext(1.25 1){$b_{m\!+\!n}$}\htext(4 1){$\cdots$}\htext(6.75
1){$b_{n\!+\!1}$} \htext(9 1){$b_n$}\htext(11.5
1){$\cdots$}\htext(14 1){$b_1$} \htext(9 3){$a_n$}\htext(11.5
3){$\cdots$}\htext(14 3){$a_1$}

\end{texdraw}}

\begin{proof} Let $M$ be a monomial in ${\mathcal M}(\la)$. Then $M$ is
expressed as
$$M=X_{a_1}(2)\cdots X_{a_n}(2)X_{b_1}(1)\cdots X_{b_{m+n}}(1)
$$ in Theorem \ref{thm:general}. We
define $\psi(M)$ by the tableau
$$\raisebox{-0.4\height}{\usebox{\tmpfiga}}\,.$$
Then it is clear that $\psi(M)$ belongs to $S(\la)$ by the
condition (i)-(iv) in Theorem \ref{thm:general}.

Conversely, for a tableau
\begin{center}
$T=\raisebox{-0.4\height}{\usebox{\tmpfiga}}\in S(\la)$,
\end{center}
we also define $\psi^{-1}(T)$ by the monomial
$$X_{a_1}(2)\cdots X_{a_n}(2)X_{b_1}(1)\cdots X_{b_{m+n}}(1).$$
Then by definition it is clear that $\psi$ and
$\psi^{-1}$ are inverses to each other.

Now, it remains to show that $\psi$ is a crystal morphism. But,
because of (\ref{eq:A1}) and (\ref{eq:A2}), it is easy to see that
$\psi$ is a crystal morphism of $U_q(G_2)$-modules from the
definition of Kashiwara operators on the set ${\mathcal M}$ of
monomials and the tensor product rule of Kashiwara operators which
is applied to the set $S(\la)$.
\end{proof}

\savebox{\tmpfiga}{\begin{texdraw} \fontsize{7}{7}\selectfont
\drawdim em \setunitscale 0.8

\rlvec(8 0) \move(0 2)\rlvec(8 0)\move(2 4)\rlvec(6 0) \move(0
0)\rlvec(0 2)\move(2 0)\rlvec(0 4)\move(4 0)\rlvec(0 4) \move(6
0)\rlvec(0 4)\move(8 0)\rlvec(0 4)

\htext(1 1){$3$}\htext(3 1){$0$}\htext(5 1){${\overline 3}$}
\htext(7 1){${\overline 2}$}

\htext(3 3){$2$}\htext(5 3){$2$}\htext(7 3){$0$}

\end{texdraw}}

\begin{example}\label{example:connection}
Let $\la$ be a dominant integral weight $\La_1+3\La_2$ and let $M$
be a monomial
$Y_1(2)^3Y_1(3)^{-3}Y_1(4)^{-1}Y_2(2)^2Y_2(3)^{-1}\in {\mathcal
M}(\la)$ given in Example \ref{example:algo} (a), then it can be
expressed as
$$
M=X_0(2)X_2(2)^2X_{\overline 2}(1)X_{\overline 3}(1)X_0(1)X_3(1).
$$
Then we have
$$\psi(M)=\raisebox{-0.4\height}{\usebox{\tmpfiga}}\,.
$$

\end{example}

We have the following proposition between $S(\la)$ and $T(\la).$

\begin{prop}\label{recovey}\, {\rm \cite{Lec}} \,\,
For a dominant integral weight $\la=m\La_1+n\La_2$, there is a
crystal isomorphism $\varphi: S(\la)\rightarrow T(\la)$ for
$U_q(G_2)$-module given by $$\varphi(S)=S_{2,1}\leftarrow
S_{2,2}\leftarrow\cdots\leftarrow S_{2,n}\leftarrow
S_{1,1}\leftarrow\cdots\leftarrow S_{1,m},$$ where $S_{i,j}\in
S(\La_i)$ $(i=1,2)$ is the column of $S$ of length $i$  from right
to left.
\end{prop}

\begin{cor} Let $\la=m\La_1+n\La_2$ be a dominant integral weight.
There is a crystal isomorphism $\phi:{\mathcal M}(\la)\rightarrow
T(\la).$
\end{cor}
\begin{proof}
By Theorem \ref{thm:general} and Proposition \ref{recovey},
$\phi=\varphi\circ \psi$ is a crystal isomorphism.
\end{proof}

\savebox{\tmpfiga}{\begin{texdraw} \fontsize{7}{7}\selectfont
\drawdim em \setunitscale 0.8

\rlvec(2 0)\rlvec(0 4)\rlvec(-2 0)\rlvec(0 -4)\move(0 2)\rlvec(2
0)

\htext(1 1){${\overline 2}$}\htext(1 3){$0$}
\end{texdraw}}
\savebox{\tmpfigb}{\begin{texdraw} \fontsize{7}{7}\selectfont
\drawdim em \setunitscale 0.8

\rlvec(2 0)\rlvec(0 4)\rlvec(-2 0)\rlvec(0 -4)\move(0 2)\rlvec(2
0)

\htext(1 1){$\overline{3}$}\htext(1 3){$2$}
\end{texdraw}}
\savebox{\tmpfigc}{\begin{texdraw} \fontsize{7}{7}\selectfont
\drawdim em \setunitscale 0.8

\rlvec(2 0)\rlvec(0 4)\rlvec(-2 0)\rlvec(0 -4)\move(0 2)\rlvec(2
0)

\htext(1 1){$0$}\htext(1 3){$2$}
\end{texdraw}}
\savebox{\tmpfigd}{\begin{texdraw} \fontsize{7}{7}\selectfont
\drawdim em \setunitscale 0.8

\rlvec(2 0)\rlvec(0 2)\rlvec(-2 0)\rlvec(0 -2)

\htext(1 1){$3$}
\end{texdraw}}
\savebox{\tmpfige}{\begin{texdraw} \fontsize{7}{7}\selectfont
\drawdim em \setunitscale 0.8

\rlvec(0 4)\rlvec(8 0)\rlvec(0 -2)\rlvec(-8 0)\move(0 0)\rlvec(6
0)\rlvec(0 4)\move(2 0)\rlvec(0 4)\move(4 0)\rlvec(0 4)

\htext(1 1){$3$}\htext(1 3){$1$}\htext(3
1){$\overline{3}$}\htext(3 3){$2$}\htext(5
1){$\overline{3}$}\htext(5 3){$3$}\htext(7 3){$\overline{1}$}
\end{texdraw}}

\begin{example}
Let $M$ be a monomial
$Y_1(2)^3Y_1(3)^{-3}Y_1(4)^{-1}Y_2(2)^2Y_2(3)^{-1}$ given in
Example \ref{example:connection}. Then we have
$$
\aligned
\phi(M)&=S_{2,1}\leftarrow S_{2,2}\leftarrow S_{2,3}\leftarrow S_{1,1}\\
&=\raisebox{-0.4\height}{\usebox{\tmpfiga}}\leftarrow
\raisebox{-0.4\height}{\usebox{\tmpfigb}}\leftarrow
\raisebox{-0.4\height}{\usebox{\tmpfigc}}\leftarrow
\raisebox{-0.3\height}{\usebox{\tmpfigd}}\\
&=\raisebox{-0.4\height}{\usebox{\tmpfige}}\,.
\endaligned
$$

\vskip 5mm Conversely, let $T$ be a tableau of $T(\La_1+3\La_2)$
$$T=\raisebox{-0.4\height}{\usebox{\tmpfige}}\,.$$
By applying the reverse bumping rule to the entries from bottom to
top and from right to left, i.e., from the entry $\overline{1}$ of
the rightmost column to the entry $1$ on top of the leftmost
column, we have the following sequence
\begin{center}
$(0,\overline{2},2,\overline{3},2,0,3)$.
\end{center}
Therefore, we have
\begin{center}
$S_{2,1}=\raisebox{-0.4\height}{\usebox{\tmpfiga}}$\,,\,\,
$S_{2,2}=\raisebox{-0.4\height}{\usebox{\tmpfigb}}$\,,\,\,
$S_{2,3}=\raisebox{-0.4\height}{\usebox{\tmpfigc}}$\, and
$S_{1,1}=\raisebox{-0.3\height}{\usebox{\tmpfigd}}$\,,
\end{center}
and since
\begin{center}
$\psi^{-1}(S_{2,1})=X_0(2)X_{\overline 2}(1)$,
$\psi^{-1}(S_{2,2})=X_2(2)X_{\overline 3}(1)$,\\
$\psi^{-1}(S_{2,3})=X_2(2)X_0(1)$, $\psi^{-1}(S_{1,1})=X_3(1)$,
\end{center}
we have
$$\aligned
\varphi^{-1}(T)&=\psi^{-1}(S_{2,1})\psi^{-1}(S_{2,2})\psi^{-1}(S_{2,3})\psi^{-1}(S_{1,1})\\
&=Y_1(2)^3Y_1(3)^{-3}Y_1(4)^{-1}Y_2(2)^2Y_2(3)^{-1}.
\endaligned$$
\end{example}

\vskip 5mm \subsection{Another realization of crystal graphs} In
this subsection, we discuss about another connected component
which is isomorphic to the crystal graph $B(\la)$ with a dominant
integral weight $\la$. More precisely, let $M_0$ be a monomial
$Y_1(-1)^{m}Y_2(-2)^{n}=X_1(-1)^{m}(X_1(-1)X_2(-2))^{n}$. Then
$M_0$ is a monomial of weight $m\La_1+n\La_2$ such that
$\eit M_0=0$ for all $i=1,2$. Moreover, we have

\begin{thm} \label{thm:another-general}
Let $\la=m\La_1+n\La_2$. Then the connected component ${\mathcal
M}'(\la)$ containing
\begin{center}
$M_0=Y_1(-1)^{m}Y_2(-2)^{n}=X_1(-1)^{m}(X_1(-1)X_2(-2))^{n}$
\end{center}
is expressed as the set of monomials
$$\begin{array}{c} M=X_{a_1}(-1)\cdots X_{a_{m+n}}(-1)
X_{b_1}(-2)\cdots X_{b_n}(-2)\\
{\rm(}a_1\preceq \cdots \preceq a_{m+n},\,\,b_1\preceq \cdots \preceq b_n{\rm)}
\end{array}$$ satisfying following conditions:
\begin{itemize}

\item [{\rm (i)}] $0$'s can not be repeated in
$a_1,\cdots,a_{m+n}$ and $b_1,\cdots,b_n$, respectively,

\item [{\rm (ii)}] $a_j\prec b_j$ unless $a_j=b_j=0$ for
$j=1,\cdots,n$,

\item [{\rm (iii)}] for $j=1,\cdots,n$,
\begin{center}
$\begin{array}{ll} {\rm dist}(a_j,b_j)\le 2 &\text{if\,
$a_j=1$ or $0$,}\\
{\rm dist}(a_j,b_j)\le 3 &\text{otherwise,}
\end{array}$
\end{center}

\item [{\rm (iv)}] for $j=1,\cdots,n-1$,
\begin{center}
$\begin{array}{ll} {\rm dist}(a_{j},b_{j+1})\ge 3 &\text{if\,
$a_j=2,3$ or $0$,}\\
{\rm dist}(a_{j},b_{j+1})\ge 2 &\text{if\, $a_j={\overline 3}$.}
\end{array}$
\end{center}
\end{itemize}
\end{thm}

\begin{proof}
By the same argument in Proposition \ref{prop:general}, it is
proved. So we omit it.
\end{proof}

\begin{thm}
Let $\la=m\La_1+n\La_2$ be a dominant integral weight. Then there
is a crystal isomorphism $\psi:{\mathcal M}'(\la)\rightarrow
T(\la)$.
\end{thm}

\savebox{\tmpfiga}{\begin{texdraw} \fontsize{7}{7}\selectfont
\drawdim em \setunitscale 0.8

\rlvec(15 0) \move(0 -2)\rlvec(15 0)\move(0 -4)\rlvec(7 0) \move(0
0)\rlvec(0 -4)\move(2 0)\rlvec(0 -4)\move(5 0)\rlvec(0 -4) \move(7
0)\rlvec(0 -4)\move(9.5 0)\rlvec(0 -2)\move(12.5 0)\rlvec(0 -2)
\move(15 0)\rlvec(0 -2)

\htext(1 -1){$a_1$}\htext(3.5 -1){$\cdots$}\htext(6
-1){$a_n$}\htext(8.25 -1){$a_{n\!+\!1}$}\htext(11
-1){$\cdots$}\htext(13.75 -1){$a_{m\!+\!n}$}\htext(1
-3){$b_1$}\htext(3.5 -3){$\cdots$}\htext(6 -3){$b_n$}
\end{texdraw}}

\begin{proof} Let $M$ be a monomial in ${\mathcal M}'(\la)$. Then $M$ is
expressed as
$$M=X_{a_1}(-1)\cdots X_{a_{m+n}}(-1)X_{b_1}(-2)\cdots X_{b_n}(-2).$$  We define $\psi(M)$ by the
tableau
$$\raisebox{-0.4\height}{\usebox{\tmpfiga}}\,.$$
Then it is clear that $\psi(M)$ belongs to $T(\la)$ by the
condition (i)-(iv) in Theorem \ref{thm:another-general}.

Conversely, for a tableau
\begin{center}
$T=\raisebox{-0.4\height}{\usebox{\tmpfiga}}\in T(\la)$,
\end{center}
we also define $\psi^{-1}(T)$ by the monomial
$$X_{a_1}(-1)\cdots X_{a_{m+n}}(-1)X_{b_1}(-2)\cdots X_{b_n}(-2).$$
Then by definition it is clear that $\psi$ and
$\psi^{-1}$ are inverses to each other. Moreover, it is also easy
to see that $\psi$ is a crystal morphism of $U_q(G_2)$-modules.
\end{proof}

\vskip 15mm


\providecommand{\bysame}{\leavevmode\hbox
to5em{\hrulefill}\thinspace}

\end{document}